\documentclass[a4paper,reqno]{amsart}
\usepackage[english]{babel}
\usepackage[T1]{fontenc} 
\usepackage[utf8]{inputenc}
\allowdisplaybreaks  % permite cortes de página en fórmulas
\usepackage{amssymb}

\usepackage{faktor}
\usepackage[svgnames,x11names,dvipsnames]{xcolor}
\usepackage{tikz-cd}
\tikzset{isom/.style={draw=none,every to/.append style={edge node={node [sloped, allow upside down, auto=false]{$\simeq$}}}}}
\usepackage{bbm}
\usepackage{mathtools}
\usepackage{microtype}
\usepackage{orcidlink}
\usepackage{hyperref}
\newif\ifdraft
\draftfalse
\usepackage{xstring}
\usepackage[abspath]{currfile}
\usepackage[yyyymmdd]{datetime}

\usepackage{enumitem}
\setlist[itemize,1]{leftmargin=*}
\setlist[enumerate,1]{leftmargin=*,label=\upshape\bfseries\arabic*.}
\setlist[enumerate,2]{leftmargin=*,align=left,label=\upshape\bfseries(\roman*),widest=iii}
\newlist{primenumerate}{enumerate}{1}
\setlist[primenumerate,1]{leftmargin=*,label=\upshape\bfseries\arabic*$'$.}
\newlist{enumerateprime}{enumerate}{1}
\setlist[enumerateprime,1]{leftmargin=*,label=\upshape\bfseries\arabic*.\phantom{$'$}}

\ifdraft
\makeatletter
\newcommand\oddfoot[1]{\gdef\@oddfoot{\reset@font#1}}
\newcommand\evenfoot[1]{\gdef\@evenfoot{\reset@font#1}}
\makeatother
\oddfoot{%
\parbox{\textwidth}{
\rule{\textwidth}{0.4pt}\par
\tiny\sffamily%
Archivo fuente: \textcolor{Red3}{\guillemotleft\StrSubstitute{\currfilename}{"}{}\guillemotright}.
Última modificación: \textcolor{Red3}{\today:\currenttime}.
}
}
\evenfoot{%
\parbox{\textwidth}{
\rule{\textwidth}{0.4pt}\par
\tiny\sffamily%
Archivo fuente: \textcolor{Red3}{\guillemotleft\StrSubstitute{\currfilename}{"}{}\guillemotright}.
Última modificación: \textcolor{Red3}{\today:\currenttime}.
}
}
\newtheoremstyle{comentario}%
{}%
{}%
{\sffamily}%
{0pt}%
{\upshape}%
{}%
{0pt}%
{\color{Red2}{\sffamily\thmname{#1}\framebox{\small\textbf{\thmnumber{#2}}}}\mdseries\upshape\thmnote{ (#3)} }%
\theoremstyle{comentario}

\theoremstyle{plain}
\RequirePackage{soulutf8}
\sethlcolor{MistyRose}
\usepackage{xparse}
\NewDocumentCommand{\badref}{v}{\textcolor{Firebrick1}{\framebox{\textbf{??}}\hl{\texttt{#1}}}}
\definecolor{labelkey}{rgb}{0.8,0.8,0.8}
\RequirePackage[outer,right,marginal]{showlabels}

\RequirePackage{xstring}
\usepackage[hyperpageref]{backref} % For drafts
\renewcommand*{\backrefalt}[4]{%
\ifcase#1%
\textcolor{Red1}{(No citations)}%
\or%
\textcolor{Yellow3}{(One citation on page #2.)}%
\else%
\textcolor{Green3}{(#1 citations in pages #2.)}%
\fi%
}
\fi

\theoremstyle{theorem}
\newtheorem{theorem}[equation]{Theorem}
\newtheorem{corollary}[equation]{Corollary}
\newtheorem{lemma}[equation]{Lemma}
\newtheorem*{lemma*}{Lemma}
\newtheorem{proposition}[equation]{Proposition}

\theoremstyle{example}
\newtheorem{definition}[equation]{Definition}
\newtheorem{remark}[equation]{Remark}
\newtheorem{example}[equation]{Example}

\numberwithin{equation}{section} %\numberwithin{defn}{section}

\DeclareMathOperator{\chr}{char}

\DeclareMathOperator{\Div}{Div}

\DeclareMathOperator{\content}{V}

\DeclarePairedDelimiterXPP\fpws[1]{}{\lbrack\mkern-1.2mu\lbrack}{\rbrack\mkern-1.2mu\rbrack}{}{\ifblank{#1}{\cdot}{#1}}
\DeclarePairedDelimiterXPP\laurent[1]{}{\lparen\mkern-3.5mu\lparen}{\rparen\mkern-3.5mu\rparen}{}{\ifblank{#1}{\cdot}{#1}}
\DeclarePairedDelimiterXPP\pairing[1]{}{\langle}{\rangle}{}{\ifblank{#1}{\cdot}{#1}}
\DeclarePairedDelimiterXPP\gen[1]{}{\langle}{\rangle}{}{\ifblank{#1}{\cdot}{#1}}
\DeclarePairedDelimiterX\braces[1]{\lbrace}{\rbrace}{\ifblank{#1}{\:\cdot\:}{#1}}
     
\newcommand{\A}{\mathbb{A}}

\newcommand{\bb}{\mathfrak{B}}

\newcommand{\f}{\mathbbm{f}}

\newcommand{\I}{\mathbb{I}}
\renewcommand{\k}{\Bbbk}

\newcommand{\m}{\mathfrak{m}}
\newcommand{\n}{\mathfrak{n}}
\newcommand{\oo}{\mathcal{O}}
\newcommand{\p}{\mathbbm{p}}
\renewcommand{\t}{\mathbbm{t}}

\let\stdphi\phi
\let\phi\varphi
\let\varphi\stdphi

\newcommand{\eqdef}{\coloneqq}

\newcommand{\suchthat}{:}
\newcommand{\AXp}[1]{\A_{X}\lbrace\mkern0mu{#1}\mkern0mu\rbrace}
\newcommand{\Kxp}[1]{K_{x}\lbrace\mkern0mu{#1}\mkern0mu\rbrace}
\newcommand{\Axp}[1]{A_{x}\lbrace\mkern0mu{#1}\mkern0mu\rbrace}
\newcommand{\AXSp}[2][S]{\A_{X,#1}\lbrace\mkern0mu{#2}\mkern0mu\rbrace}

\newcommand{\joinrelshort}{\mathrel{\mkern-8mu}}
\newcommand{\shortlongrightarrow}{\relbar\joinrelshort\rightarrow}
\newcommand{\isomto}{\mathrel{\mathop{\setbox0\hbox{$\mathsurround0pt\shortlongrightarrow$}\ht0=0.3\ht0\box0}\limits^{\hspace{-1pt}\scalebox{1.0}{$\sim$}\mkern2mu}}}

\newcommand*{\prodprime}{\operatornamewithlimits{%
\mathchoice
{\prod\nolimits\raisebox{1.618ex}{\hspace{0em}\makebox[0pt]{$\scriptstyle\prime$}}\hspace{-0.2em}}%display
{\prod\nolimits\raisebox{1ex}{\hspace{-0.15em}\makebox{$\scriptstyle\prime$}}\hspace{0.1em}}%text
{\prod\nolimits\raisebox{0.618ex}{\hspace{-0.16em}\makebox{$\scriptstyle\prime$}}\hspace{0.1em}}%script
{\prod\nolimits\raisebox{0.3ex}{\hspace{-0.2em}\makebox{$\scriptstyle\prime$}}\hspace{0em}}%scriptscript
}}

\title[Characterization of subfields of adelic algebras]{%
Characterization of subfields of adelic algebras by a product formula
}%
\author[L.~M.~Navas~Vicente]{Luis Manuel Navas Vicente\orcidlink{0000-0002-5742-8679}}
\email{navas@usal.es}

\author[F.~J.~PLaza~Mart\'in]{Francisco J. Plaza Mart\'in\orcidlink{0000-0002-5532-7567}}
\email{fplaza@usal.es}

\thanks{Research of both authors supported by grant PID2023-150787NB-I00 and of the first author also by grant PID2021-124332NB-C22 (Gobierno de España). ORCID ID: 0000-0002-5742-8679, 0000-0002-5532-7567.}

\thanks{The authors have no competing interests to declare. No datasets are associated to this research.}
\address{Departamento de Matem\'aticas and IUFFYM, Universidad de
Salamanca,  Plaza de la Merced 1-4
        \\
        37008 Salamanca. Spain.
        \\
         Tel: +34 923294460. 
}

\begin{document}
%%%%%%%%%%%%%%%%%%%%%%%%%%%%%%%%%%%%%%%%%%%%%%%%%%%%%%%%%%%%%%%%%%%%%%%%%%%%%%%%%%%
% OPCIONES LATEX GENERALES
%%%%%%%%%%%%%%%%%%%%%%%%%%%%%%%%%%%%%%%%%%%%%%%%%%%%%%%%%%%%%%%%%%%%%%%%%%%%%%%%%%%
\emergencystretch 3em 
%%%%%%%%%%%%%%%%%%%%%%%%%%%%%%%%%%%%%%%%%%%%%%%%%%%%%%%%%%%%%%%%%%%%%%%%%%%%%%%%%%%
%%%%%%%%%%%%%%%%%%%%%%%%%%%%%%%%%%%%%%%%%%%%%%%%%%%%%%%%%%%%%%%%%%%%%%%%%%%%%%%%%%%

\begin{abstract}
We consider projective, irreducible, non-singular curves over an algebraically closed field $\k$. A cover $Y \to X$ of such curves corresponds to an extension $\Omega/\Sigma$ of their function fields and yields an isomorphism $\A_{Y} \simeq \A_{X} \otimes_{\Sigma} \Omega$ of their geometric adele rings. The primitive element theorem shows that $\A_{Y}$ is a quotient of $\A_{X}[T]$ by a polynomial. 

In general, we may look at quotient algebras $\AXp{\p} = \A_{X}[T]/(\p(T))$ where $\p(T) \in \A_{X}[T]$ is monic and separable over $\A_{X}$, and try to characterize the field extensions $\Omega/\Sigma$ lying in $\AXp{\p}$ which arise from covers as above. We achieve this topologically, namely, as those $\Omega$ which embed discretely in $\AXp{\p}$, and in terms of an additive analog of the product formula for global fields, a result which is reminiscent of classical work of Artin-Whaples and Iwasawa.

The technical machinery requires studying which topology on $\AXp{\p}$ is natural for this problem. Local compactness no longer holds, but instead we have linear topologies defined by commensurability of $\k$-subspaces which coincide with the restricted direct product topology with respect to integral closures. The content function is given as an index measuring the discrepancy in commensurable subspaces.
\ifdraft	
\leavevmode
\begin{description}

\item[\sffamily Archivo fuente] {\sffamily\textcolor{Red3}{\currfilename}}

\item[\sffamily Última modificación] {\sffamily\textcolor{Red3}{\today:\currenttime}.}

\end{description}
\fi
\end{abstract}

\maketitle

% \footnotesize
% \tableofcontents
% \normalsize

%%%%%%%%%%%%%%%%%%%%%%%%%%%%%%%%%%%%%%%%%%%%%%%%%%%%%%%%%%%%%%%%%%%%%%%%%%%%%%%%%%%
%%%%%%%%%%%%%%%%%%%%%%%%%%%%%%%%%%%%%%%%%%%%%%%%%%%%%%%%%%%%%%%%%%%%%%%%%%%%%%%%%%%
%%%%%%%%%%%%%%%%%%%%%%%%%%%%%%%%%%%%%%%%%%%%%%%%%%%%%%%%%%%%%%%%%%%%%%%%%%%%%%%%%%%
\section{Introduction}
\label{sec:introduction}
% (fold)
%%%%%%%%%%%%%%%%%%%%%%%%%%%%%%%%%%%%%%%%%%%%%%%%%%%%%%%%%%%%%%%%%%%%%%%%%%%%%%%%%%%
%%%%%%%%%%%%%%%%%%%%%%%%%%%%%%%%%%%%%%%%%%%%%%%%%%%%%%%%%%%%%%%%%%%%%%%%%%%%%%%%%%%
%%%%%%%%%%%%%%%%%%%%%%%%%%%%%%%%%%%%%%%%%%%%%%%%%%%%%%%%%%%%%%%%%%%%%%%%%%%%%%%%%%%

In~\cite{MPNPP-collectanea} we studied some aspects of the ring $\A_{X}$ of geometric adeles of an algebraic curve $X$ over an algebraically closed base field $\k$. In particular, we discussed the topological structure, in analogy with the well-known classical case of adeles of global fields. A fundamental difficulty is that the topology is no longer locally compact, hence we lose the machinery of Haar measure and ultrametric harmonic analysis as introduced by Tate. Instead there is a linear topology characterized by commensurability of $\k$-subspaces (the original idea of commensurability is also due to Tate in~\cite{TateRes}).
Finiteness of the cohomology implies that the function field $\Sigma$ of the curve still embeds discretely in the adele ring as in the classical case, although other well-known results no longer hold. Nevertheless, we were able to establish results regarding reciprocity laws on the curve, connecting with some of the themes introduced in~\cite{AndersonPablos,MP,MP-mediterranean}.

Let $\Sigma$ be the function field of $X$. An irreducible cover $\pi : Y \to X$ corresponds to a field extension $\Omega/\Sigma$. The adele ring $\A_{Y}$ is then topologically isomorphic to $\A_{X} \otimes_{\Sigma} \Omega$ as an $\A_{X}$-algebra (this is analogous to the case of global fields discussed in~\cite{CaFr}). Together with the primitive element theorem for fields, this implies that $\A_{Y}$ is a quotient of $\A_{X}[T]$ by a polynomial.

% The inverse problem consists in determining when an $\A_{X}$-algebra of the form $\AXp{\p} = \A_{X}[T]/(\p(T))$ where $\p(T)$ is a monic separable polynomial with coefficients in $\A_{X}$, is isomorphic to the adele ring $\A_{Y}$ arising from a cover. We have solved this problem for Kummer extensions in~\cite{adeles03}.

This raises the following question. Consider an $\A_{X}$-algebra of the form $\AXp{\p} = \A_{X}[T]/(\p(T))$ where $\p(T)$ is a monic separable polynomial with coefficients in $\A_{X}$. These algebras can also be endowed with a topological structure, characterized in terms of commensurability of $\k$-subspaces, in a manner suggested by~\cite{MP} and generalizing the case of $\A_{X}$. Now, we may look for subfields $\Omega$ embedded in $\AXp{\p}$ which come from covers as above. Namely, given a finite extension of fields $\Omega/\Sigma$ with $\Omega$ inside $\AXp{\p}$, and $Y$ the Zariski-Riemann variety of $\Omega$, can we determine when the geometric adele ring $\A_{Y}$ embeds in $\AXp{\p}$ or is in fact isomorphic to it. 

Our answer to this question is contained in the following main results. Theorem~\ref{T:equivalences inyectivity} gives a topological characterization stating that the natural map $\A_{Y} \to \AXp{\p}$ is an embedding if and only if $\Omega$ is discrete. Next, Theorem~\ref{T:product formula equivalence} shows that an equivalent condition is that $\Omega$ satisfy a ``product formula'', which in this case is actually an additive version. This is expressed via an additive content function we introduce on the group of invertible elements $\AXp{\p}^{*}$, deriving from a notion of index measuring the discrepancy between commensurable subspaces. This is reminiscent of the classical results of~\cite{ArtinW}, where global fields are axiomatically characterized by the existence of a product formula for valuations, and also of~\cite{Iwa}, where a topological characterization of 
the adele ring of a global field is given.

% (end)
% end section introduction

%%%%%%%%%%%%%%%%%%%%%%%%%%%%%%%%%%%%%%%%%%%%%%%%%%%%%%%%%%%%%%%%%%%%%%%%%%%%%%%%%%%
%%%%%%%%%%%%%%%%%%%%%%%%%%%%%%%%%%%%%%%%%%%%%%%%%%%%%%%%%%%%%%%%%%%%%%%%%%%%%%%%%%%
%%%%%%%%%%%%%%%%%%%%%%%%%%%%%%%%%%%%%%%%%%%%%%%%%%%%%%%%%%%%%%%%%%%%%%%%%%%%%%%%%%%
\section{Adelic algebras}
\label{sec:adelic algebras}
%%%%%%%%%%%%%%%%%%%%%%%%%%%%%%%%%%%%%%%%%%%%%%%%%%%%%%%%%%%%%%%%%%%%%%%%%%%%%%%%%%%
%%%%%%%%%%%%%%%%%%%%%%%%%%%%%%%%%%%%%%%%%%%%%%%%%%%%%%%%%%%%%%%%%%%%%%%%%%%%%%%%%%%
%%%%%%%%%%%%%%%%%%%%%%%%%%%%%%%%%%%%%%%%%%%%%%%%%%%%%%%%%%%%%%%%%%%%%%%%%%%%%%%%%%%

%%%%%%%%%%%%%%%%%%%%%%%%%%%%%%%%%%%%%%%%%%%%%%%%%%%%%%%%%%%%%%%%%%%%%%%%%%%%%%%%%%%
%%%%%%%%%%%%%%%%%%%%%%%%%%%%%%%%%%%%%%%%%%%%%%%%%%%%%%%%%%%%%%%%%%%%%%%%%%%%%%%%%%%
%%%%%%%%%%%%%%%%%%%%%%%%%%%%%%%%%%%%%%%%%%%%%%%%%%%%%%%%%%%%%%%%%%%%%%%%%%%%%%%%%%%
\subsection{Geometric adeles}
\label{subsec:adeles}
% (fold)
%%%%%%%%%%%%%%%%%%%%%%%%%%%%%%%%%%%%%%%%%%%%%%%%%%%%%%%%%%%%%%%%%%%%%%%%%%%%%%%%%%%
%%%%%%%%%%%%%%%%%%%%%%%%%%%%%%%%%%%%%%%%%%%%%%%%%%%%%%%%%%%%%%%%%%%%%%%%%%%%%%%%%%%
%%%%%%%%%%%%%%%%%%%%%%%%%%%%%%%%%%%%%%%%%%%%%%%%%%%%%%%%%%%%%%%%%%%%%%%%%%%%%%%%%%%

We start by giving a brief summary of the algebraic and topological construction of the adeles of a curve. Although the construction is analogous to the case of global fields given in~\cite{CaFr}, there are some differences in the geometric case, especially regarding the topology. See~\cite[\S2]{MPNPP-collectanea} for more details about geometric adeles and~\cite{Warner} for the general theory of topological rings. The following discussion is a brief summary of the former.

Let $X$ be a projective, irreducible, non-singular curve over an algebraically closed field $\k$. Let $\Sigma$ be the function field of $X$. We fix the following notation:
\begin{itemize}

% \item $\Sigma$ is the function field of $X$.

\item When we write $x \in X$ it will be implicitly assumed that $x$ is a closed point of $X$, corresponding to the valuation $\upsilon_{x}$ on the function field $\Sigma$. Denote by $\oo_{X,x}$, or simply $\oo_{x}$, the valuation ring at $x$. Since $\k$ is algebraically closed, the closed points are in one-to-one correspondence with the discrete valuations on $\Sigma$.

\item Let $A_{x}$ be the completion of $\oo_{x}$ with respect to $\upsilon_{x}$, which will also denote the extended valuation. Let $\m_{x}$ be the maximal ideal of $A_{x}$, $K_{x}$ its quotient field. Since $\k$ is algebraically closed, the residue field is $A_{x}/\m_{x} = \k$. 

\end{itemize}
%%
% A choice of uniformizing parameter $z_{x}$ at $x$ determines the following isomorphisms:% (\cite{Cohen}):
% %%
% \begin{equation}
% \label{E:Ax Kx pws}
% 	A_{x} \simeq \k\fpws{z_{x}},
% 	\quad
% 	K_{x} \simeq \k\laurent{z_{x}},
% 	\quad
% 	A_{x}^{*} \simeq \k^{*} \times (1 + \m_{x}).
% \end{equation}
% %%

The ring of \emph{adeles} $\A_{X}$ of $\Sigma/\k$ is the subring of $\prod_{x} K_{x}$ given by the restricted direct product with respect to the subrings $A_{x}$, over the closed points of $X$,
\begin{equation*}
\label{E:adele}
	       \A_{X} 
	\eqdef \prodprime_{x \in X} (K_{x},A_{x})
	 = \left\{ (\alpha_{x})_{x\in X} \suchthat \alpha_{x}\in A_{x}  \mbox{ for almost all } x\in X \right\}
\end{equation*}
where ``almost all'' means ``for all but finitely many''. It is equipped with the restricted product topology, which we recall is \emph{not} induced by the direct product topology.

% For an adele $\alpha = (\alpha_{x})_{x\in X}$ and a (closed) point $x \in X$ with $\alpha_{x}\in A_{x}$, $\alpha(x)$ will denote the image of its $x$-component $\alpha_{x}$ in the residue field $\k(x)$.

$\A_{X}$ arises as a direct limit as follows: denote by $S$ a finite subset of (closed) points of $X$. Consider
\begin{equation}
\label{E:AXS}
		   \A_{X,S}
	\eqdef \prod_{x \in S} K_{x}
	\times \prod_{x \in X \setminus S} A_{x}.
\end{equation}
This ring carries a linear topology; namely, the one generated by the neighborhood basis at $0$ consisting of the sets $\prod_{x} \m_{x}^{n_{x}}$, where $(n_{x})_{x\in X}$ runs over collections of non-negative integers $n_{x}$ such that $n_{x}=0$ for almost all $x$.

Given finite subsets $S_1,S_2$ with $S_1\subseteq S_2$, we have an inclusion $\A_{X,S_1}\hookrightarrow \A_{X,S_2}$. Fix a finite subset $S_{0}$ and consider sets containing $S_{0}$. Then the direct limit with respect to these inclusions is isomorphic to $\A_{X}$:
\begin{equation}
\label{E:AX limit}
	\A_{X} \simeq \varinjlim_{S \supseteq S_{0}} \A_{X,S},
\end{equation}
as linearly topological rings. In particular, this does not depend on the choice of $S_{0}$.

\begin{remark}
\label{R:direct limits and restricted direct product}
The identification of the direct limit with the restricted direct product holds in the most general framework where these objects are defined, including the case of topological groups, rings, etc. We will make use of this fact several times.
\end{remark}

There is an alternative way to generate this topology in the adele ring. Let $\A_{X}^+$ denote the subring of $\A_{X}$ given by the usual direct product
\[
	\A_{X}^+ \eqdef \prod_{x \in X} A_{x}, 
\]
i.e. where the word ``almost'' is dropped in the previous definition. Recalling the general notion of integrality for ring extensions, it is easy to check that $\A_{X}^{+}$ is integrally closed in $\A_{X}$, highlighting the interaction of the algebraic notion of integral closure in commutative rings, with the adelic topology. We will come back to this theme in \S\ref{subsec:algebraic+topological structure AXp}. Consider the following notion.

\begin{definition}
\label{D:commensurability}
Two vector subspaces $U,V$ of a vector space over a field $\k$ are called \emph{commensurable}, denoted by $U\sim V$, iff 
\[
	\dim_{\k} (U+V) / (U\cap V)  <  \infty.
\]
It is easy to see that commensurability is an equivalence relation on $\k$-subspaces.
\end{definition}

One may check that the topology of $\A_{X}$ coincides with the topology having  $\{U\subseteq \A_{X}   \suchthat   U \sim \A_{X}^+\}$ as neighborhood basis at $0$. It is worth noticing that $\A_{X}^{+} \subseteq \A_{X,S}$ for all $S$.

Finally, the \emph{idele group} $\I_{X}$ is the group $\A_{X}^{*}$ of invertible elements of $\A_{X}$ endowed with the initial topology of the map $\I_{X}  \rightarrow   \A_{X} \times \A_{X}$ that sends $\alpha$ to $(\alpha, \alpha^{-1})$. Observe that it is the restricted product of $K_{x}^{*}$ with respect to the unit groups $A_{x}^{*}$ and that $\I_X$ can be constructed as the direct limit $\varinjlim \I_{X,S}$ where $\I_{X,S} \eqdef \A^{*}_{X,S}$.

\begin{remark}
\label{R:Sigma discrete in AX}
From now on, we shall always consider the geometric adeles $\A_{X}$ as a topological ring with the topology described above. An important fact to keep in mind is that the function field $\Sigma$ embedded in $\A_{X}$ then has the discrete topology, just as is the case for global fields. This follows from noting that if $V \subseteq \A_{X}$ is a $\k$-vector space, then $V$ is discrete if and only if $\dim_{\k} (V \cap \A_{X}^{+}) < \infty$ (in which case $V$ is automatically closed). For $V = \Sigma$ the finite-dimensionality follows from observing that $\Sigma \cap \A_{X}^{+} = H^{0}(X,\oo_{X}) = \k$ under our hypotheses on $X$. 
\end{remark}

For a more detailed exposition of the above results, especially regarding the equivalence of the adelic (restricted direct product) and the topology of commensurability, see~\cite[\S2.2]{MPNPP-collectanea}.

% (end)
% end subsection adeles

%%%%%%%%%%%%%%%%%%%%%%%%%%%%%%%%%%%%%%%%%%%%%%%%%%%%%%%%%%%%%%%%%%%%%%%%%%%%%%%%%%%
%%%%%%%%%%%%%%%%%%%%%%%%%%%%%%%%%%%%%%%%%%%%%%%%%%%%%%%%%%%%%%%%%%%%%%%%%%%%%%%%%%%
\subsection{Definitions and notation}
\label{subsec:definition of adelic algebra}
%%%%%%%%%%%%%%%%%%%%%%%%%%%%%%%%%%%%%%%%%%%%%%%%%%%%%%%%%%%%%%%%%%%%%%%%%%%%%%%%%%%
%%%%%%%%%%%%%%%%%%%%%%%%%%%%%%%%%%%%%%%%%%%%%%%%%%%%%%%%%%%%%%%%%%%%%%%%%%%%%%%%%%%

% For the remainder of the paper, $n$ will denote a fixed integer prime to $\chr \k$. When $n$ is assumed prime, it will be denoted by $p$.

Consider an irreducible cover $Y \to X$ or, equivalently, a finite separable field extension $\Omega$ of the function field $\Sigma$ of $X$. The existence of primitive elements implies that there is a monic separable polynomial $P(T) \in \Sigma[T]$ such that $\Omega \simeq \Sigma[T]/(P(T))$. Since $\Sigma$ embeds in the adele ring $\A_{X}$, it can be shown that tensoring yields
\begin{equation}
\label{E:cassels isom quotient algebraic}
		\faktor{\A_{X}[T]}{(P(T))} \simeq \A_{X} \otimes_{\Sigma} \Omega
	                  \simeq \A_{Y},
\end{equation}
where the second isomorphism is analogous to~\cite[(14.2)]{CaFr} for global fields. In other words, the adele ring of a cover $Y$ is a polynomial quotient algebra over $\A_{X}$ of the above form. This provides the motivation for much of what follows.

It is important to recall that the classical result cited above is actually a topological isomorphism. Thus one of the things we need to do is introduce topologies
on $\A_{X}$-algebras in the geometric case. This is done in the next section. Before turning to this, we introduce notation and deal with the important algebraic aspect of separability.

% An important class of algebras of the form $\AXp{\p}$ arises as the adele rings of  irreducible covers $Y \to X$, or equivalently, finite separable field extensions $\Omega$ of the function field $\Sigma$ of $X$. Indeed, the existence of primitive elements implies that there is a monic separable polynomial $\p(T) \in \Sigma[T] \hookrightarrow \A_{X}[T]$ such that $\Omega \simeq \Sigma\{\p\} = \Sigma[T]/(\p(T))$. Arguing as in~\cite[(14.2)]{CaFr} for global fields, ignoring the topological structure for the moment, we have
% \[
% 	\faktor{\A_{X}[T]}{(\p(T))} \simeq \A_{X} \otimes_{\Sigma} \Omega
% 	                  \simeq \A_{Y}.
% \]
% The converse problem of when a given algebra $\AXp{\p}$ is of this form is taken up in \S\ref{sec:subfields of adelic algebras}.

In general, given a commutative ring $R$, an indeterminate $T$, and a polynomial $P(T) \in R[T]$, we will denote the corresponding polynomial quotient algebra by
\begin{equation}
\label{E:Rxp}
	R\{P\} \eqdef \faktor{R[T]}{(P(T))}.
\end{equation}
We are focused on the case $R = \A_{X}$, which we distinguish notationally by double letters. Thus we will study $\A_{X}$-algebras of the form
\[
	\AXp{\p} = \faktor{\A_{X}[T]}{(\p(T))}
	\suchthat
	\p(T) \in \A_{X}[T].
\]
% As we have seen, an important class of these are the adele rings of irreducible covers of $X$.
One of our main goals is to answer the question of when a given algebra of the form $\AXp{\p}$ is the adele ring of a cover $Y \to X$. We take this up in \S\ref{sec:subfields of adelic algebras} after introducing the various objects and structures we need.

To a polynomial $\p(T) \in \A_{X}[T]$ with adelic coefficients and a point $x \in X$, it is natural to associate the polynomial $\p_{x}(T) \in K_{x}[T]$ whose coefficients are the projections onto $K_{x}$ of the coefficients of $\p(T)$. For almost all $x$, in fact $\p_{x}(T) \in A_{x}[T]$. For these reasons we are led to consider also
% Thus we will also need to look at
% In particular, we are interested in the following three cases: $R = \A_{X}$ and $R = K_{x}$ or $A_{x}$ for some $x \in X$. We will distinguish them notationally by using double letters for adeles. Thus we will consider
%
\begin{itemize}

% \item $\AXp{\p} = \A_{X}[T] / (\p(T))$, where $\p(T) \in \A_{X}[T]$.

\item $\Kxp{P} = K_{x}[T] / (P(T))$, for $x \in X$ and $P(T) \in K_{x}[T]$.

\item $\Axp{P} = A_{x}[T] / (P(T))$, for $x \in X$ and $P(T) \in A_{x}[T]$.

\end{itemize}
%
% \begin{definition}
% \label{D:algebra AXp}
% Let $\p(T) \in \A_{X}[T]$, where $T$ is an indeterminate, be a polynomial with coefficients in $\A_{X}$. Consider the $\A_{X}$-algebra
% %%
% \begin{equation}
% \label{E:AXp}
% 	\AXp{\p} \eqdef \A_{X}[T] / (\p(T)).
% \end{equation}
% %%
% \end{definition}
% 
% 

In~\cite{adeles02} and~\cite{adeles03} we studied the special class of polynomials of the form $\p(T) = T^{n} - \t$ where $\t \in \I_{X}$ is an idele (to ensure separability). Algebraically this leads to the Kummer theory of $\A_{X}$, and geometrically to the study of cyclic covers of $X$.

In view of the underlying geometric motivation and also to simplify the technical details, we will assume that $\p(T)$ is monic of some fixed degree $n \geq 1$ and separable over $\A_{X}$. In particular this guarantees that each $\p_{x}$ is also monic of degree $n$, independently of $x$.

We need to make use of the general theory of separable algebras over commutative rings, keeping in mind that $\A_{X}$ is a ``bad'' ring, insofar as it has infinitely many idempotents and zero divisors, so care must be taken when using results (e.g.~\cite{Ford, Jan, Nagahara}) on irreducibility and separability, since many criteria cannot be directly applied to adele rings.

Given a commutative ring $R$, a commutative $R$-algebra $A$ is said to be separable over $R$ if $A$ is projective as an $A\otimes_{R} A$-module. If $P(T)\in R[T]$ is a monic polynomial, we say that $P(T)$ is separable over $R$ if $R[T]/(P(T))$ is a separable $R$-algebra. We have the following characterization of separable polynomials.

\begin{proposition}[{\cite[Prop. 4.6.1]{Ford}}]
\label{P:caracterizationseparable}
A monic polynomial $P(T) \in R[T]$ over a commutative ring $R$ is separable if and only if the ideal generated by $P$ and its formal derivative $P'$ contains $1$, i.e. $(P(T),P'(T))=R[T]$.
	%\item $\Delta(P)\in R^*$ where $\Delta(P)$ denotes the discriminant of the polynomial $P(T)$.
\end{proposition}

% \begin{proof}
% See \cite[Proposition 4.6.1]{Ford}.
% %The hypothesis show that $S\eqdef  R[T]/(P(T))$ is an $R$-module free of rank equal to $m\eqdef  \operatorname{deg}P$. The equivalence of $(1)$ and $(2)$ is given in \cite[Proposition 4.6.1]{Ford}. Note that the discriminant $\Delta(P)$ equals the determinant of the trace form of the extension $S$ taken with respect to the basis $\{1,T, \ldots, T^{m-1}\}$ and, hence, it generates $D_R(S)$ The equivalence of $(1)$ and $(3)$ follows from \cite[Thm. 4.7]{Knus}
% \end{proof}

% The previous result is also a consequence of a theorem of Nagahara~\cite[Theorem 2.3]{Nagahara}, which can also be applied to relate the separability of a polynomial over $\A_{X}$ with the separability of its projections over the residue fields of $\A_{X}$. However, because of the rather intricate structure of the maximal spectrum of $\A_{X}$ (see~\cite{SerranoHolgado}), we aim at a characterization only in terms of the closed points of $X$, or equivalently, in terms of the projections at each point $x \in X$, which, by definition, have all their coefficients in $A_{x}$ at almost every $x$.

For a monic polynomial with coefficients in $\A_{X}$, we may relate the separability of $\p(T) \in \A_{X}[T]$ to that of its projections $\p_{x}(T) \in K_{x}[T]$.

\begin{definition}\label{D:pointwiseseparability}
A monic polynomial $\p(T)\in\A_{X}[T]$ is said to be \emph{pointwise separable} if  $\p_{x}(T)$ is separable over $K_{x}$ for every $x\in X$.
\end{definition}

\begin{proposition}\label{P:sep-point-sep}
A monic polynomial $\p(T)\in\A_{X}[T]$ is separable if and only if is pointwise separable and for almost every $x \in X$, $\p_{x}(T) \in A_{x}[T]$ is separable over $A_{x}$.
\end{proposition}

\begin{proof}
%Given $P(T)$, let $\Delta(P(T))$ denote its discriminant.  Having in mind Proposition~\ref{P:caracterizationseparable}, is follows that $\Delta(P(T))\in \A_{X}^* = \I_X$ if and only if $\Delta(P(T)_{x})= \Delta(P(T))_{x} \in K_{x}^*$ for all $x\in X$  and $\Delta(P(T)_{x})  \in A_{x}^*$ for almost all $x\in X$.
This is a matter of noting that a Bézout identity $\mathbbm{a}(T)\p(T)+\mathbbm{b}(T)\p'(T)=1$ in $\A_{X}[T]$ implies a corresponding Bézout identity $\mathbbm{a}_{x}(T) \p_{x}(T) + \mathbbm{b}_{x}(T) \p'_{x}(T) = 1$ in $K_{x}[T]$ at each $x \in X$, which is actually an identity in $A_{x}[T]$ at almost every $x$.
Conversely, given such identities at each $x$, since $K_{x}$ is the field of fractions of $A_{x}$, they can be chosen so that $\deg \mathbbm{a}_{x}(T)$ and $\deg \mathbbm{b}_{x}(T)$ are bounded by $\deg \p_{x}(T)$, which is independent of $x$, equal to $\deg \p(T)$ since $\p(T)$ was chosen to be monic. Hence the pointwise Bézout identities can be combined to give an identity in $\A_{X}[T]$.
\end{proof}

\begin{remark}
\label{R:points classification}
In light of Proposition~\ref{P:sep-point-sep}, we see that for a monic separable polynomial $\p(T) \in \A_{X}[T]$ and a point $x \in X$, one of the following three cases will hold: either $\p_{x}(T)$ does has some coefficient in $K_{x}$ but not $A_{x}$, or $\p_{x}(T) \in A_{x}[T]$ but it is not separable over $A_{x}$, or $\p_{x}(T) \in A_{x}[T]$ and it is separable over $A_{x}$. The latter holds for almost all points. We will return to this point in~Proposition~\ref{P:Hensel complete splitting}.
\end{remark}

\begin{example}
\label{E:separability Tn-t}
$\p(T) = T^{n} - \t$ where $\t \in \A_{X}$, is separable if and only if in fact $\t \in \I_{X}$ (\cite[Proposition 2.10]{adeles02}).
\end{example}

\subsection{Algebraic and topological structure of adelic algebras}
\label{subsec:algebraic+topological structure AXp}
% (fold)
%%%%%%%%%%%%%%%%%%%%%%%%%%%%%%%%%%%%%%%%%%%%%%%%%%%%%%%%%%%%%%%%%%%%%%%%%%%%%%%%%%%
%%%%%%%%%%%%%%%%%%%%%%%%%%%%%%%%%%%%%%%%%%%%%%%%%%%%%%%%%%%%%%%%%%%%%%%%%%%%%%%%%%%

In order to study the structure of $\A_{X}$-algebras of the form $\AXp{\p}$, where $\p(T) \in \A_{X}[T]$ is monic and separable, we introduce topologies analogous to those discussed in~\S\ref{subsec:adeles}. The resulting topological $\A_{X}$-algebra structure is highly relevant for characterizing the subfields of $\AXp{\p}$ in \S\ref{sec:subfields of adelic algebras}.

%%%%%%%%%%%%%%%%%%%%%%%%%%%%%%%%%%%%%%%%%%%%%%%%%%%%%%%%%%%%%%%%%%%%%%%%%%%%%%%%%%%
\subsubsection{Adelic algebras as restricted direct products}
\label{subsubsec:AXp as restricted direct product}
% (fold)
%%%%%%%%%%%%%%%%%%%%%%%%%%%%%%%%%%%%%%%%%%%%%%%%%%%%%%%%%%%%%%%%%%%%%%%%%%%%%%%%%%%

From now on we will consider a separable monic polynomial $\p(T) \in \A_{X}[T]$ and the associated separable adelic algebra $\AXp{\p} = \A_{X}[T]/(\p(T))$. Let us see how it can be expressed as a restricted direct product. Note that
\begin{equation}
\label{E:Axp tensor Kx}
		\Kxp{\p_{x}} = \faktor{K_{x}[T]}{(\p_{x}(T))} = \AXp{\p} \otimes_{A_{x}} K_{x}
\end{equation}
since $K_{x}$ is the field of fractions of $A_{x}$. In addition, since the coefficients of $\p(T)$ are adeles, for almost all $x$, $\p_{x}(T) \in A_{x}[T]$, thus
\[
	\Axp{\p_{x}} = \faktor{A_{x}[T]}{(\p_{x}(T))} 
\]
is well-defined and can be regarded as a subring of $\Kxp{\p_{x}}$.

Since $S_0 = \{x \in X \suchthat \p_{x}(T) \notin A_{x}[T]\}$ is finite, we may define, for a finite subset $S \supseteq S_{0}$ of closed points,
\begin{equation}
\label{E:AXSp}
			   \AXSp{\p}
	\eqdef \prod_{x \in S} \Kxp{\p_{x}}
	\times \prod_{x \in X \setminus S} \Axp{\p_{x}}
	\simeq \A_{X,S}[T]/(\p(T)).
\end{equation}

\begin{proposition}
\label{P:AXp restricted product Axpx}
The adelic algebra $\AXp{\p}$ is the restricted direct product of $\Kxp{\p_{x}}$ with respect to the subrings $\Axp{\p_{x}}$, i.e.
\begin{equation}
\label{E:AXp restricted product Axpx}
    \AXp{\p} \simeq \varinjlim_{S \supseteq S_{0}} \AXSp{\p}
	         = \prodprime_{x \in X} (\Kxp{\p_{x}}, \Axp{\p_{x}}).
\end{equation}
\end{proposition}

\begin{proof}
Take the direct limit in~\eqref{E:AXSp} with respect to the obvious inclusion maps between the $\AXSp{\p}$.
\end{proof}

An alternative description of $\AXp{\p}$ as a restricted direct product may be given which is more convenient for several reasons. It consists of considering, instead of $\AXp{\p_{x}}$, the integral closure of $A_{x}$ in $\Kxp{\p_{x}}$. This will be denoted by $\bb_{x}$. One advantage is that $\bb_{x}$ is defined for all $x$, obviating the need to take into account the exceptional set $S_{0}$ above. In addition, in the geometric case of an algebra arising from a cover $\pi : Y \to X$, if $\pi(y) = x$, the completion of the local ring $\oo_{Y,y}$ is precisely the integral closure of $\pi^{*} A_{x}$. This also holds for the adeles of a global field.

To verify this claim we need the following result, which is interesting in itself and will be used several times.

\begin{proposition}
\label{P:Hensel complete splitting}
Let $P(T) \in A_{x}[T]$ be monic of degree $n$. Then the following conditions are equivalent: 
\begin{enumerate}

\item $P(T)$ is separable over $A_{x}$.

\item The discriminant $\Delta = \prod_{i<j} (a_{i} - a_{j}) \in A_{x}$ is invertible.

\item The reduced polynomial $\overline{P}(T) \in \k[T]$ is separable.

\end{enumerate}

\noindent If this is the case, then:

\begin{enumerate}[resume]

\item  $P(T)$ splits completely over $A_{x}$, i.e. $P(T) = (T-a_{1}) \cdots (T-a_{n})$ with $a_{1}, \dots, a_{n} \in A_{x}$.

\item The isomorphism
\begin{equation}
\label{E:isomorphism to copies of Kx}
\begin{split}
	&\psi_{x} : \Kxp{P} = \faktor{K_{x}[T]}{(P(T))} \simeq \prod^{n} K_{x}
	\\
	&\psi_{x}(Q(T)) = (Q(a_{1}),Q(a_{2}),\dots,Q(a_{n}))
\end{split}
\end{equation}
induces an isomorphism of $\Axp{P}$ with $\prod^{n} A_{x}$.

\item The integral closure $\bb_{x}$ of $A_{x}$ in $\Kxp{P}$ is $\Axp{P}$.

\end{enumerate}
\end{proposition}

\begin{proof}
The equivalence of the first three conditions follows from~\cite[Theorem~2.3]{Nagahara} and the fact that the invertible elements of $A_{x}$ are those with zero valuation. Assuming that this holds, we then have:
\begin{enumerate}

\setcounter{enumi}{3}

\item This follows from Hensel's Lemma, taking into account that $A_{x}$ is a complete DVR with algebraically closed residue field $\k$, and $(P(T),P'(T)) = A_{x}[T]$ by Proposition~\ref{P:caracterizationseparable}.

\item It is clear that $\psi_{x}$ is an isomorphism and embeds $\Axp{P}$ in $\prod^{n} A_{x}$, thus it suffices to check that it is surjective. This is a well-known standard argument: given $(b_0, \dots ,b_{n-1}) \in \prod^{n} A_{x}$, we need to find a polynomial $\sum c_{i} T^{i} \in A_{x}[T]$ of degree less than $n$ such that 
\[
	c_{0} + c_{1} a_{k} + \dots + c_{n-1} a_{k}^{n-1} = b_{k}
	\quad
	(k = 0,1,\dots,n-1).
\]
This is a system of equations in the $c_{j}$, with (Vandermonde) determinant $\Delta$, therefore by (2) there is in fact a solution with $c_{i} \in A_{x}$.

\item Simply note that integral closure commutes with finite products and $A_{x}$-algebra isomorphisms, $A_{x}$ is integrally closed in $K_{x}$, and thus $\bb_{x} = \psi_{x}^{-1}(\prod^{n} A_{x}) = \Axp{P}$.
\qedhere
\end{enumerate}
\end{proof}

% The result of the lemma is also valid without the assumption that $e_{x} = 1$, replacing $A_{x}$ with the ring of integers of the cyclic extension $K_{x}(\tau_{x}^{1/e_{x}})$ where $\tau_{x} = t_{x}^{1/m_{x}}$ as in Lemma~\ref{L:structureKxntx}. We do not need this greater generality.

% Note that, for the time being, we make no mention of topological structure. This will be done later in \S\ref{subsec:algebraic+topological structure AXp}. The previous remarks lead us also naturally to consider the algebraic structure of $\Kxp{\p_{x}}$ and $\Axp{\p_{x}}$ as $x$ varies over the (closed) points of $X$.

% To define a restricted direct product topology on $\AXp{\p}$, in view of its definition and Proposition~\ref{P:AXp restricted product Axpx}, the logical choice of subrings would seem to be the $\Axp{\p_{x}}$.

\begin{corollary}
\label{C:AXp direct limit+restricted product}
Given a monic separable polynomial $\p(T) \in \A_{X}[T]$, the adelic algebra $\AXp{\p}$ is isomorphic to the restricted direct product of the rings $\Kxp{\p_{x}}$ with respect to the integral closures $\bb_{x}$ of $A_{x}$ in $\Kxp{\p_{x}}$, that is to say,
\begin{equation}
\label{E:AXp restricted product Kxp Bx}
           \AXp{\p} 
	\simeq \varinjlim_{S} \biggl(\prod_{x \in S} \; \Kxp{\p_{x}} \times \!\! \prod_{x \in X \setminus S} \bb_{x} \biggr)
    = \prodprime_{x \in X} (\Kxp{\p_{x}}, \bb_{x})
\end{equation}
where $S$ runs over all finite subsets of $X$.
\end{corollary}

\begin{proof}
By Proposition~\ref{P:AXp restricted product Axpx}, it is enough to show that in fact $\Axp{\p_{x}} = \bb_{x}$ for almost all $x$. This follows from Proposition~\ref{P:sep-point-sep} and Proposition~\ref{P:Hensel complete splitting}.
\end{proof}

Corollary~\ref{C:AXp direct limit+restricted product} allows us to express an element $\alpha \in \AXp{\p}$ in the form $\alpha = (\alpha_{x})_{x \in X}$ with $\alpha_{x} \in \Kxp{\p_{x}}$ and $\alpha_{x} \in \bb_{x}$ for almost all $x$.

% (end)
% end subsubsection AXp as restricted direct product

%%%%%%%%%%%%%%%%%%%%%%%%%%%%%%%%%%%%%%%%%%%%%%%%%%%%%%%%%%%%%%%%%%%%%%%%%%%%%%%%%%%
\subsubsection{The restricted direct product topology}
\label{subsubsec:the restricted direct product topology}
% (fold)
%%%%%%%%%%%%%%%%%%%%%%%%%%%%%%%%%%%%%%%%%%%%%%%%%%%%%%%%%%%%%%%%%%%%%%%%%%%%%%%%%%%

We now turn our attention to defining a natural topology on adelic algebras of the form $\AXp{\p}$. This necessitates defining topologies on each factor $\Kxp{\p_{x}}$ which define a restricted direct product topology. It will coincide with the corresponding induced direct limit topology via the identification in~\eqref{E:AXp restricted product Kxp Bx} (see Remark~\ref{R:direct limits and restricted direct product}). In addition we will see that there is an equivalent characterization via the notion of commensurability for $\k$-vector subspaces.

\begin{lemma}
\label{L:nx-adic topology}
There is a unique Hausdorff topology on $\bb_{x}$ (resp. $\Kxp{\p_{x}}$) restricting to the $\m_{x}$-adic topology on $A_{x}$ (resp. $K_{x}$). It is the $\n_{x}$-adic topology, where $\n_{x}$ is the ideal of $\bb_{x}$ equal to the integral closure of $\m_{x}$ in $\Kxp{\p_{x}}$.
\end{lemma}

\begin{proof}
Given any $x \in X$, recall (Proposition~\ref{P:sep-point-sep}) that $\p_{x}$ is separable over $K_{x}$. Thus $\Kxp{\p_{x}}$ is an étale $K_{x}$-algebra, therefore isomorphic to a finite product of field extensions of $K_{x}$. Noting that taking integral closure commutes with finite direct products and $K_{x}$-algebra isomorphisms, the statement reduces to the case where $\Kxp{\p_{x}}$ is a finite field extension $L$ of $K_{x}$, which is well-known.
\end{proof}

% %%
% \begin{equation}
% \label{E:Kxp product Kxj}
% 	\Kxp{\p_{x}} \simeq \prod_{j=1}^{m_{x}} K_{x,j}
% \end{equation}
% %%
% for some $m_{x} \geq 1$. In addition, for $j=1,\dots,m_{x}$, let $B_{x,j}$ be the integral closure of $A_{x}$ in $K_{x,j}$. Note that each  $B_{x,j}$ is a DVR, and denote its maximal ideal by $\m_{x,j}$. We then have
% %%
% \begin{equation}
% \label{E:Bx decomposition Bxj}
% 		\bb_{x} = \psi_{x}^{-1}(B_{x,1} \times \dots \times B_{x,m_{x}}).
% \end{equation}
% %%

We will henceforth assume given the $\n_{x}$-adic topology on each $\Kxp{\p_{x}}$. Observe that this in fact gives $\Kxp{\p_{x}}$ the structure of a topological $K_{x}$-algebra with respect to which the integral closure $\bb_{x}$ is open.

\begin{proposition}
\label{P:adelic mx-adic topology on AXp}
The $\n_{x}$-adic topologies defined in Lemma~\ref{L:nx-adic topology} induce a linear topology on $\AXp{\p}$ with respect to which it is a topological $\A_{X}$-algebra. A neighborhood basis at $0$ for the topology on $\AXp{\p}$ is given by the collection of sets
\begin{equation}
\label{E:UD neighborhood}
		U_{D} \eqdef \prod_{x\in X}\n_{x}^{k_{x}},
\end{equation}
where $D = \sum k_{x} x$ is an effective divisor on $X$. We will refer to this topology on $\AXp{\p}$ simply as the adelic topology.
\end{proposition}

\begin{proof}
In view of Corollary~\ref{C:AXp direct limit+restricted product}, a collection of topologies on each $\Kxp{\p_{x}}$ with respect to which the subrings $\bb_{x}$ are open, determines a restricted direct product topology on $\AXp{\p}$, with a neighborhood basis at $0$ given by the sets $U = \prod_{x \in S} U_{x} \times \prod_{x \notin S} \bb_{x}$, where $S \subseteq X$ is a finite subset and $U_{x}$ is a neighborhood of $0$ in $\Kxp{\p_{x}}$. By Lemma~\ref{L:nx-adic topology}, we may assume without loss of generality that $U_{x} = \n_{x}^{k_{x}}$ for some $k_{x} \geq 0$, and thus $U = U_{D}$ for $D = \sum_{x} k_{x} x$.

Note that for an effective divisor $D = \sum k_{x} x$ on $X$, by definition we have $\A_{X} \cap U_{D} = \prod_{x \in X} \m_{x}^{k_{x}}$, and the latter sets are a neighborhood basis at $0$ for the usual adelic topology on $\A_{X}$. This shows that the action of $\A_{X}$ on $\AXp{\p}$ is continuous with respect to their corresponding adelic topologies.
\end{proof}

It is useful to note that the neighborhoods $U_{D}$, for $D \geq 0$, satisfy:
\begin{itemize}
%
% \item $U_{D'} = \cup_{D\geq D'} U_{D}$.
%
% \item $U_{D} \subseteq U_{D'}\cap U_{D''}$ for any $D$ with $D\geq D'$, $D\geq D''$ and, in particular, $U_{D}\subseteq U_{D'}$ for  $D\geq D'$,
%
\item $\cap_{D} U_{D} = (0)$.

\item $U_{D} \subseteq U_{D'}\cap U_{D''} \subseteq U_{D'} + U_{D''}$ for $D, D', D'' $ with $D\geq D'$, $D\geq D''$.

% \item In particular if $D\geq D'$,  then $U_{D}\subseteq U_{D'}$ and $\dim_{k}(U_{D'}/U_{D})< \infty$.
%
\end{itemize}
In addition, the $U_{D}$ are actually $\k$-vector subspaces of $\AXp{\p}$.
We denote the neighborhood corresponding to the divisor $D=0$ by
\begin{equation}
\label{E:AXplus}
		\AXp{\p}^+ \eqdef \prod_{x\in X} \bb_{x}.
\end{equation}

The following observation is immediate from the proof, but we point it out since a generalization forms the basis of a characterization of discrete subfields of $\AXp{\p}$ (see Theorem~\ref{T:equivalences inyectivity}).

\begin{corollary}
\label{C:adelic topology on AX in AXp is initial topology}
The adelic topology on $\A_{X}$ is the initial topology induced by the inclusion $\A_{X} \hookrightarrow \AXp{\p}$. In particular, the inclusion is continuous.
\end{corollary}

In a similar manner as mentioned above for $\A_{X}$, there is a further useful characterization of the direct product topology on $\AXp{\p}$ in terms of the notion of commensurability (Definition~\ref{D:commensurability}).

\begin{lemma}
\label{L:conmensurability topology equivalence}
The set of $\k$-vector subspaces $U$ of $\AXp{\p}$ commensurable with $\AXp{\p}^{+}$ is a neighborhood basis at $0$ for the restricted direct product topology.
\end{lemma}

\begin{proof}
We may use the same reasoning as in~\cite[Prop. 2.1]{MPNPP-collectanea}, namely, it is easy to check that for $D \geq 0$, the subspace $U_{D}$ is commensurable with $\AXp{\p}^{+} = U_{0}$, i.e., that $\dim_{\k}(U_{D}+U_{0})/(U_{D}\cap U_{0})<\infty$. 

Conversely, if $U$ is a subspace with $\dim_{\k}(U+\AXp{\p}^{+})/(U\cap \AXp{\p}^{+})<\infty$ then, since $(0) = \bigcap_{D} U_{D}$, there is some divisor $D\in \Div(X)$ such that $U_{D}\subset U$, thus $U$ is a neighborhood of $0$ in the restricted direct product topology.
\end{proof}

% (end)
% end subsubsection the restricted direct product topology

%%%%%%%%%%%%%%%%%%%%%%%%%%%%%%%%%%%%%%%%%%%%%%%%%%%%%%%%%%%%%%%%%%%%%%%%%%%%%%%%%%%
\subsubsection{Global structure}
\label{subsubsec:global structure}
% (fold)
%%%%%%%%%%%%%%%%%%%%%%%%%%%%%%%%%%%%%%%%%%%%%%%%%%%%%%%%%%%%%%%%%%%%%%%%%%%%%%%%%%%

\begin{proposition}
\label{P:AXp+ is integral closure of AX+ in AXp}
The subring $\AXp{\p}^+$ is the integral closure of $\A_{X}^+$ in $\AXp{\p}$.
\end{proposition}

\begin{proof}
The inclusion of the integral closure in $\AXp{\p}^+$ is trivial. For the reverse inclusion, note that at every point $x$, any element $\alpha_{x}$ of $\bb_{x}$ satisfies a monic polynomial with coefficients in $A_{x}$. By Proposition~\ref{P:Hensel complete splitting}, at almost every point, $\bb_{x}$ is a free $A_{x}$-module of rank $n$ and we may choose this polynomial to be the characteristic polynomial of $\alpha_{x}$ over $A_{x}$, which has fixed degree $n$. Thus for $\alpha = (\alpha_{x})_{x} \in \AXp{\p}^+$ we can find polynomials for each $\alpha_{x}$ with uniformly bounded degrees, which implies that $\alpha$ satisfies a monic polynomial with coefficients in $\A_{X}^{+}$.
\end{proof}

\begin{corollary}
\label{C:algebraic isomorphism AXT is topological}
Any $\A_{X}$-algebra isomorphism between $\AXp{\p_{1}}$ and $\AXp{\p_{2}}$ is also a topological isomorphism.
\end{corollary}

\begin{proof}
This follows from Proposition~\ref{P:AXp+ is integral closure of AX+ in AXp} and Lemma~\ref{L:conmensurability topology equivalence}. Indeed, if $\phi:\AXp{\p_{1}}\rightarrow\AXp{\p_{2}}$ is an $\A_{X}$-algebra isomorphism, it is enough to show that $\phi(\AXp{\p_{1}}^+)$ and $\AXp{\p_{2}}^+$ are commensurable. In fact, we have that $\phi(\AXp{\p_{1}}^+)= \AXp{\p_{2}}^+$, since $\phi$ preserves the integral closure of $\A_{X}^+$.
\end{proof}

\begin{theorem}\label{T:Axp is free AX-module}
Let $\p(T) \in \A_{X}[T]$ be monic and separable. Then $\A_{X}^{\oplus n} \simeq \AXp{\p} $ as  linear topological $\A_{X}$-modules. 
\end{theorem}

We start with a following technical lemma of the ``AKLB'' type.

\begin{lemma}
\label{L:commensurabilityDVR}
Let $A$ be a complete DVR with residue field $k$ and quotient field $K$. Assume furthermore that $k$ embeds in $A$. Let $L$ be an étale $K$-algebra and $B$ the integral closure of $A$ in $L$. If $y_1,\dots, y_n$ is a $K$-vector space basis of $L$, then $\dim_{k} (B+\sum A y_i )/(B\cap \sum A y_i )<\infty$.
\end{lemma}

\begin{proof}
For each $i =1,\dots,n$, since $y_i$ is algebraic over $K$, there is some $a_i\in A$ such that $y'_i \eqdef  a_i y_i$ is integral over $A$ and therefore $y'_i\in B$. Then $\sum A y'_i \subseteq B$ and $B /(\sum A y'_i )$ is a torsion $A$-module. Since $B$ is finitely generated, both sides of the exact sequence 
\[
	B /(\sum A y'_i ) \longrightarrow 
	(B+\sum A y_i )/(B\cap \sum A y_i ) \longrightarrow
	(B+\sum A y_i )/B \to 0 
\]
are finitely generated torsion $A$-modules and thus, by the classification theorem for modules over PIDs and the hypotheses on $A$, they are finite dimensional $k$-vector spaces. Hence the central term is also finite dimensional, which proves the claim.
\end{proof}

\begin{proof}[Proof of Theorem~\ref{T:Axp is free AX-module}]
Since $\p(T)$ is assumed monic, associating to a class $A(T) \mod \p(T)$ the remainder of $A(T)$ under division by $\p(T)$ gives an $\A_{X}$-module isomorphism $\psi : \Axp{\p} \isomto \A_{X}^{\oplus n}$.

% and separable, for almost all $x \in X$, the projection $\p_{x}(T)$ has coefficients in $A_{x}$ and is monic and separable over $A_{x}$. Thus Proposition~\ref{P:Hensel complete splitting} applies, yielding a $K_{x}$-algebra isomorphism $\psi_{x} : \Kxp{\p_{x}} \isomto \prod^n K_{x}$ sending $\Axp{\p_{x}}$, which in this case coincides with the integral closure $\bb_{x}$, to $\prod^{n} A_{x}$.
%
% At the finite number of exceptional points (i.e. where the hypotheses of Proposition~\ref{P:Hensel complete splitting} fail for $\p_{x}$), since $\Kxp{\p_{x}}$ is an étale $K_{x}$-algebra, and thus isomorphic to a product of extension fields of $K_{x}$, there will be many $K_{x}$-linear isomorphisms (not necessarily algebra morphisms) $\psi_{x} : \Kxp{\p_{x}} \isomto \prod^n K_{x}$, and we simply choose one. Thinking of all the $\psi_{x}$ as $A_{x}$-module isomorphisms, observe that the collection $\{\psi_{x}\}_{x \in X}$ induces a map on the restricted direct product of the $(\Kxp{\p_{x}},\bb_{x})$ and thus defines a global $\A_{X}$-module isomorphism $\psi : \AXp{\p} \isomto \A_{X}^{\oplus n}$.

We need to check that $\psi$ is a homeomorphism. In view of Lemma~\ref{L:conmensurability topology equivalence}, it is easy to check that this is equivalent to showing that $\mathfrak{A} \eqdef  \psi^{-1}((\A_{X}^{+})^{\oplus n})$ is commensurable with $\mathfrak{B} \eqdef  \AXp{\p}^{+}$, i.e. 
\[
	\dim_{\k} (\mathfrak{A}+\mathfrak{B})/(\mathfrak{A} \cap \mathfrak{B})<\infty.
	\tag{$\ast$}
\]	
Since $\p(T)$ is also assumed separable, by Proposition~\ref{P:sep-point-sep} and Proposition~\ref{P:Hensel complete splitting}, for almost all points $x$, the integral closure of $\bb_{x}$ of $A_{x}$ in $\Kxp{\p_{x}}$ is $\Axp{\p_{x}}$. Thus~$(\ast)$ follows by applying Lemma~\ref{L:commensurabilityDVR} at the remaining finite number of exceptional points $x$, with $A = A_{x}$, $K= K_{x}$, $L = \Kxp{\p_{x}}$ and $B = \bb_{x}$ and $\{y_{i} = \psi_{x}^{-1}(T^{i})\}_{i=1}^{n}$ where $\psi_{x}$ is the induced map $\psi_{x} : \Kxp{\p_{x}} \isomto K_{x}^{\oplus n}$.
\end{proof}
\section{Subfields of adelic algebras}
\label{sec:subfields of adelic algebras}
% (fold)
%%%%%%%%%%%%%%%%%%%%%%%%%%%%%%%%%%%%%%%%%%%%%%%%%%%%%%%%%%%%%%%%%%%%%%%%%%%%%%%%%%%
%%%%%%%%%%%%%%%%%%%%%%%%%%%%%%%%%%%%%%%%%%%%%%%%%%%%%%%%%%%%%%%%%%%%%%%%%%%%%%%%%%%
%%%%%%%%%%%%%%%%%%%%%%%%%%%%%%%%%%%%%%%%%%%%%%%%%%%%%%%%%%%%%%%%%%%%%%%%%%%%%%%%%%%

As we mentioned before, the prototypical example of an adelic algebra of the form $\AXp{\p}$ introduced in \S\ref{subsec:definition of adelic algebra} is the geometric adele ring of a cover of algebraic curves. We now turn to the reverse problem; namely, given an adelic algebra $\AXp{\p}$, we wish to determine whether there are separable field extensions $\Omega$ of $\Sigma$ embedded in $\AXp{\p}$ such that $\A_{X} \otimes_{\Sigma} \Omega$ embeds in $\AXp{\p}$. As we will see, to achieve this goal, it is necessary to incorporate the topological structure.

The main results are Theorem~\ref{T:equivalences inyectivity}, which characterizes such subfields $\Omega$, and Theorem~\ref{T:product formula equivalence} which relates this question to an analog of the Artin-Whaples product formula as in~\cite{ArtinW}.

%%%%%%%%%%%%%%%%%%%%%%%%%%%%%%%%%%%%%%%%%%%%%%%%%%%%%%%%%%%%%%%%%%%%%%%%%%%%%%%%%%%
%%%%%%%%%%%%%%%%%%%%%%%%%%%%%%%%%%%%%%%%%%%%%%%%%%%%%%%%%%%%%%%%%%%%%%%%%%%%%%%%%%%
%%%%%%%%%%%%%%%%%%%%%%%%%%%%%%%%%%%%%%%%%%%%%%%%%%%%%%%%%%%%%%%%%%%%%%%%%%%%%%%%%%%
\subsection{Discrete subfields}
\label{subsec:discrete subfields}
% (fold)
%%%%%%%%%%%%%%%%%%%%%%%%%%%%%%%%%%%%%%%%%%%%%%%%%%%%%%%%%%%%%%%%%%%%%%%%%%%%%%%%%%%
%%%%%%%%%%%%%%%%%%%%%%%%%%%%%%%%%%%%%%%%%%%%%%%%%%%%%%%%%%%%%%%%%%%%%%%%%%%%%%%%%%%
%%%%%%%%%%%%%%%%%%%%%%%%%%%%%%%%%%%%%%%%%%%%%%%%%%%%%%%%%%%%%%%%%%%%%%%%%%%%%%%%%%%

The first step is to show, as we mentioned in \S\ref{subsec:definition of adelic algebra}, that if we start from a cover $Y \to X$, we obtain topological isomorphisms of $\A_{X}$-algebras as expressed in~\eqref{E:cassels isom quotient algebraic}, namely
\[
		\faktor{\A_{X}[T]}{(\p(T))} \simeq \A_{X} \otimes_{\Sigma} \Omega
	                  \simeq \A_{Y},
\]
where $\p(T) \in \A_{X}[T]$ is monic and separable. Here, $\A_{Y}$ has the usual adelic topology as in \S\ref{subsec:adeles}, the middle term has the tensor product topology, and $\AXp{\p}$ has the topology described in \S\ref{subsec:algebraic+topological structure AXp}.

A cover $Y \to X$, where $Y$ is a projective, irreducible, non-singular curve, determines a finite separable field extension $\Omega$ of $\Sigma$. Conversely, given such an $\Omega$, the Zariski-Riemann variety $Y$ of $\Omega$, whose set of closed points is the set of discrete valuations on $\Omega$, determines a cover of $X$.

Given $y\in Y$, we use the following notation:
\begin{itemize}

\item $L_y$, the completion of $\Omega$ with respect to $\upsilon_y$, the valuation corresponding to $y$.

\item $B_{y}$, the completion of the local ring $\oo_{Y,y}$, with respect to $\upsilon_y$.

\item $\pi:Y\to X$, the map corresponding to $\Sigma\hookrightarrow \Omega$ which maps a valuation $\upsilon$ of $\Omega$ to its restriction to $\Sigma$.

\item The rest of the notation in Section~\ref{subsec:adeles} also applies: thus $\A_{Y}$ is the adele ring of $Y$, $\A_{Y}^{+} = \prod_{y} B_{y}$, etc.

\item In this geometric situation, $\p_{\Omega}(T) \in \Sigma[T] \subseteq \A_{X}[T]$ is the minimal polynomial of a primitive element of $\Omega/\Sigma$. Thus it is monic, separable and irreducible over $\Sigma$.

\end{itemize}
%%

% In other words, the adelic algebras of the form $\AXp{\p}$ arising from a cover $Y \to X$ as above, are algebraically and topologically isomorphic to the geometric adele ring $\A_{Y}$.

\begin{theorem}
\label{T:cassels}
Let $\Omega$ be a finite separable field extension of $\Sigma$ and let $\A_{Y}$ denote the ring of adeles of $\Omega$. Then
\begin{equation}
\label{E:cassels}
	\AXp{\p_{\Omega}} \simeq \A_{X}\otimes_{\Sigma} \Omega  \simeq   \A_{Y}	
\end{equation}
as linear topological $\A_{X}$-algebras. In particular, $\Omega$ embeds discretely in $\AXp{\p_{\Omega}}$. Furthermore, $\A_{Y}^{+}$ is the integral closure of $\A_{X}^{+}$ in $\A_{X}\otimes_{\Sigma} \Omega$.
\end{theorem}

\begin{proof}
The topology on the tensor product is assumed to be the product topology, keeping in mind that $\dim_{\Sigma} \Omega < \infty$ (this is described in~\cite[p.54]{CaFr}).	
It is straightforward to verify from the various characterizations of the restricted direct product topology, that the first isomorphism is topological.

For the second isomorphism,	one checks that the proof of the analogous statement for global fields given in~\cite[(14.2)]{CaFr} carries over to the geometric case; namely, taking restricted direct products, considered as direct limits (which commute with tensor products) in the relations
\begin{equation}
\label{E:tensorCompletion}
		K_{x} \otimes_{\Sigma} \Omega \simeq \prod_{\pi(y)=x} L_{y},
\end{equation}
yields~\eqref{E:cassels} as an algebraic isomorphism. 

Thus, to check that this second isomorphism is also topological is equivalent to showing that it is continuous and open. For this it suffices to see that the subring generated by $\A_{X}^{+}$ in $\A_{X} \otimes_{\Sigma} \Omega \simeq \A_{Y}$ is commensurable with $\A_{Y}^{+}$. Note that at any point, the integral closure of $A_{x}$ in $\prod_{\pi(y)=x} L_{y}$ is $\prod_{\pi(y)=x} B_{y}$. At an unramified point of the cover $Y \to X$, this coincides with the subring of $K_{x} \otimes_{\Sigma} \Omega$ generated by $A_{x}$. Hence it suffices to check commensurability at the remaining finitely many ramified points. For these,  Lemma~\ref{L:commensurabilityDVR} applies. The previous fact about the integral closure of $A_{x}$ also implies that the integral closure of $\A_{X}^{+}$ in $\A_{Y}$ is $\A_{Y}^{+}$.

Finally, the fact that $\Omega$ embeds discretely in $\AXp{\p_{\Omega}}$ follows immediately from the previous topological isomorphisms and Remark~\ref{R:Sigma discrete in AX}.
\end{proof}

Theorem~\ref{T:cassels} is the analog of~\cite[Lemma 8]{Iwa}, where linear local compactness is here replaced by commensurability, and shows that discreteness is a basic necessary condition for the adele ring of a field extension $\Omega/\Sigma$ to embed in some adelic algebra $\AXp{\p}$. We now consider the converse problem, fixing an adelic algebra $\AXp{\p}$ beforehand and considering varying $\Omega$.

% In the case of a Galois cover $Y \to X$, that is to say, with $\Omega/\Sigma$ a finite Galois field extension, by considering~\eqref{E:tensorCompletion} at each fiber, we can easily verify that the isomorphism $\A_{Y} \simeq \A_{X} \otimes_{\Sigma} \Omega$ from~\eqref{E:cassels} is compatible with the action of $\Gal(Y/X)$, and therefore
% %%
% \begin{equation}
% \label{E:AY Gal(Y/X) invariant is AX}
% 	\A_{Y}^{\Gal(Y/X)} = \A_{X}.
% \end{equation}
% %%

% In \S\badref{subsec:the main example}, we have seen how an adelic algebra arises from a finite separable field extension $\Omega$ of $\Sigma$ (recall~\eqref{E:cassels} and Example~\badref{EX:rational Kummer Sigma}).

% We use the notation introduced in \S\ref{subsec:definition of adelic algebra} and refer to it for the topological structure of adelic algebras.

Recalling the results and notation of \S\ref{subsec:algebraic+topological structure AXp} (especially Lemma~\ref{L:nx-adic topology}), for a given $x \in X$, the étale $K_{x}$-algebra $\Kxp{\p_{x}}$ decomposes into a product of a finite number $m_{x} \geq 1$ of finite separable field extensions of $K_{x}$. We label these by $K_{x,j}$, for $1 \leq j \leq m_{x}$. Each $K_{x,j}$ is a complete DVR, with valuation denoted by $\upsilon_{x,j}$. Denote its ring of integers by $B_{x,j}$, and its maximal ideal my $\m_{x,j}$. Thus
\begin{equation}
\label{E:KxKj BxBj}
	\Kxp{\p_{x}} \simeq \prod_{j=1}^{m_{x}} K_{x,j},
	\qquad
	\bb_{x} = \prod_{j=1}^{m_{x}} B_{x,j}	
\end{equation}
since integral closure commutes with finite direct products.

Suppose now that $\Omega$ be a finite separable field extension of $\Sigma$ embedded in $\AXp{\p}$, and $Y$ is its Zariski-Riemann variety. Let us introduce some notation. Given $x \in X$ and $1 \leq j \leq m_{x}$,~\eqref{E:Kx tensor Omega -> Kxtn} and~\eqref{E:KxKj BxBj} give us a composite map
\[
	\prod_{\pi(y) = x} L_{y} \to \prod_{i=1}^{m_{x}} K_{x,i} \to K_{x,j}
\]
whose kernel is a prime ideal in a product of fields, and thus maximal, equal to a product of $L_{y'}$ with $y' \in \pi^{-1}(x) \setminus \{y\}$ for a unique $y \in \pi^{-1}(x)$. In other words, the map factors through a unique $L_{y}$. We will denote this point by
\[
	y \eqdef s_{j}(x).
\]
% By definition, the valuation $\upsilon_{x,j}$ on $K_{x,j}$ restricts to $\upsilon_{y}$ on both $\Omega$ and $L_{y}$.
% 

% various $\upsilon_{x,j}$ corresponding to different indices $j$ can restrict to the same valuation on $\Omega$.

% For each $j : 1 \leq j \leq m_{x}$, the composition
% %%
% \begin{equation}
% \label{E:Omega into AXp into Kxj}
% 		\Omega \hookrightarrow \AXp{\p} \to K_{x,j}
% \end{equation}
% %%
% of the inclusion and the projection induces an inclusion of fields, and hence a valuation on $\Omega$ by restriction. The uniqueness of extensions of valuations in DVRs justifies using the same notation $\upsilon_{x,j}$ for it.

There is a natural map from the adele ring $\A_{Y}$, which as we have seen is isomorphic to $\AXp{\p_{\Omega}}$, to $\AXp{\p}$, which we proceed to study.

\begin{lemma}
\label{L:AY -> AXp continuity}
Let $\Omega \subset \AXp{\p}$ be a finite separable field extension of $\Sigma$ and  $Y$ is its Zariski-Riemann variety. The inclusion $\Omega \hookrightarrow \AXp{\p}$ factors through a continuous map
\begin{equation}
\label{E:AY -> AXp}
	\A_{Y} \xrightarrow{\phi} \AXp{\p}.
\end{equation}
\end{lemma}

\begin{proof}
Consider the inclusions $\Sigma\hookrightarrow \Omega\hookrightarrow \AXp{\p}$. Completing at a point $x \in X$ and taking tensor products, we obtain maps
\begin{equation}
\label{E:Kx tensor Omega -> Kxtn}
	       K_{x} \otimes_{\Sigma} \Omega 
	\simeq K_{x} \otimes_{\Sigma} \A_{Y}
	\simeq \prod_{\pi(y)=x} L_{y}
	\xrightarrow{\phi_{x}} K_{x} \otimes_{\A_{X}} \AXp{\p}
	=      \Kxp{\p_{x}}.	
\end{equation}
The integral closure of $A_{x}$ in $\prod_{\pi(y)=x} L_{y}$ is $\prod_{\pi(y)=x} B_{y}$. The latter map sends this product to the integral closure $\bb_{x}$ of $A_{x}$ in $\Kxp{\p_{x}}$. Taking restricted direct products and recalling that $\A_{X} \otimes_{\Sigma} \Omega \simeq \A_{Y}$, the maps $\phi_{x}$ yield the desired map $\phi : \A_{Y} \to \AXp{\p}$. Continuity follows from the fact that the inverse image of $\AXp{\p}^{+}$ is open in $\A_{Y}$, which is straightforward.
\end{proof}

\begin{lemma}
\label{L:AY -> AXp injectivity}
The map $\phi$ in~\eqref{E:AY -> AXp} is injective if and only if
\begin{equation}
\label{E:AY -> AXp injectivity}
	j \mapsto s_{j}(x) : \{ 1 \leq j \leq m_{x} \} \to \pi^{-1}(x)
\end{equation}
is surjective, for all $x \in X$.
\end{lemma}

\begin{proof}
Observe that~\eqref{E:AY -> AXp} is injective if and only if the map~\eqref{E:Kx tensor Omega -> Kxtn} is injective for each $x \in X$ which, in light of~\eqref{E:KxKj BxBj} and the definition of $s_{j}(x)$, happens if and only if for every $y\in\pi^{-1}(x)$ there is at least one $j : 1 \leq j \leq m_{x}$ such that $y = s_{j}(x)$. 
\end{proof}

Note that, even if $\phi$ is injective, it is certainly possible that $s_{j}(x) = s_{j'}(x)$ for different indices $j,j'$, e.g., in the trivial case $\Omega = \Sigma$, $s_{j}(x) = x$ for all $j$,  while the factorization~\eqref{E:KxKj BxBj} depends on that of $\p_{x}(T)$ over $K_{x}$ (recall for example the case discussed in Proposition~\ref{P:Hensel complete splitting}).

% ~\eqref{E:AXp restricted product Kxp Bx}

% Corolarios 5.9, 5.14 de Atiyah-MacDonald.

\begin{theorem}[Characterization of discrete subfields]
\label{T:equivalences inyectivity}
Let $\p(T) \in \A_{X}[T]$ be monic and separable over $\A_{X}$. For a given finite separable field extension $\Omega$ of $\Sigma$ contained in  $\AXp{\p} = \A_{X}[T]/(\p(T))$, the following statements are equivalent:
\begin{enumerate}

\item $\Omega$ is discrete in $\AXp{\p}$.

\item The map~\eqref{E:AY -> AXp} is injective: $\A_{Y} \hookrightarrow \AXp{\p}$.

\item The adelic topology of $\A_{Y}$ coincides with the initial topology induced by the map~\eqref{E:AY -> AXp}.

\end{enumerate}
If these conditions hold, then $\deg_{\Sigma}\Omega\leq n$.
\end{theorem}

\begin{proof}
\leavevmode
\begin{itemize}

\item $(1) \implies (2)$:
If $\A_{Y} \to \AXp{\p}$ is not injective, by Lemma~\ref{L:AY -> AXp injectivity}, there is some $y\in Y$ such that if $x = \pi(y)$, then $y \neq s_{j}(x)$ for all $j = 1,\dots,m_{x}$.

Now, let $D= \sum k_{x} x$ be an effective divisor on $X$, and $U_D= \prod\n_{x}^{k_{x}}$ be the neighbourhood of $0$ in $\AXp{\p}$ it defines. 
Here, $\n_{x}$ is the ideal defined in Lemma~\ref{L:nx-adic topology}. From the proof of this result it follows that $\n_{x}$ is the ideal of $\bb_{x}$ corresponding under~\eqref{E:KxKj BxBj} to the product $\m_{x,1} \times \dots \times \m_{x,m_{x}}$ of all the maximal ideals.

By the Riemann-Roch Theorem, there is a nonconstant rational function $f_{D} \in \Omega$ which is regular in $Y\setminus\{y\}$ and lies in $U_D$; that is to say, $\upsilon_{y'}(f_{D}) \geq 0$ for all $y' \neq y$ and $\upsilon_{s_{j}(x)}(f_{D}) \geq k_{x}$ for each $(x,j)$.

Since this holds for every effective divisor $D$ on $X$, we have that $\Omega\cap U_D\neq\{0\}$ for every $D \geq 0$, therefore $\Omega$ is not discrete in $\AXp{\p}$.

% \item $(2) \implies (1)$: Suppose that $\A_{Y}\to\AXp{\p}$ is injective. We have to show that there is an effective nonzero divisor $D= \sum k_{x} x$ in $X$ such that $\Omega\cap U_D= \{0\}$. Otherwise, since for every $y\in Y$ there is some $(x,j)$ such that $y=s_j(x)$, we would obtain a nonzero function on $\Omega$ with zeroes and no poles, which is impossible.

\item $(2) \implies (3)$: Given $x \in X$, by~\eqref{E:Kx tensor Omega -> Kxtn}, we have the commutative square
\begin{equation}
\label{E:injectivity AY AXp}
\begin{tikzcd}
	\A_{Y} \arrow[r,hook,"\phi"] \arrow[d] & \AXp{\p} \arrow[d]
	\\
	\displaystyle\prod_{\pi(y)=x} L_{y} \arrow[r,hook,"\phi_{x}"'] & \Kxp{\p_{x}}.
\end{tikzcd}
\end{equation}
The integral closure of $A_{x}$ in $\prod_{\pi(y)=x} L_{y}$ is  $\prod_{\pi(y)=x} B_{y}$. Via~\eqref{E:KxKj BxBj}, we may identify the integral closure $\bb_{x}$ of $A_{x}$ in $\Kxp{\p_{x}}$ with $\prod_{j=1}^{m_{x}} B_{x,j}$; thus
\begin{equation}
	\begin{aligned}
	\prod_{\pi(y)=x} B_{y}
	&\overset{\eqref{E:AY -> AXp injectivity}}{=} \phi_{x}^{-1}\biggl(\prod_{j=1}^{m_{x}} B_{x,j}\biggr)
	\\
	&= \phi_{x}^{-1}\biggl(\bb_{x} \cap \phi_{x}\Bigl(\prod_{\pi(y)=x} L_{y}\Bigr)\biggr)
	= \phi_{x}^{-1}(\bb_{x}) \cap \prod_{\pi(y)=x} L_{y}
	\end{aligned}
\label{E:injectivity Bx prod Ly = prod By}
\end{equation}
for all $x \in X$. Thus recalling~\eqref{E:AXplus}, we have
\[
	\begin{aligned}
	\A_{Y}^{+} &= \prod_{y \in Y} B_{y}
	           =  \prod_{x \in X} \prod_{\pi(y) = x} B_{y}
			    = \phi^{-1}\biggl(\prodprime_{x} \, \bb_{x} \cap \phi_{x}\Bigl(\prod_{\pi(y)=x} L_{y}\Bigr)\biggr)
			  \\
	          &= \phi^{-1}(\AXp{\p}^{+} \cap \phi(\A_{Y}))
			   = \phi^{-1}(\AXp{\p}^{+}) \cap \A_{Y}.
	\end{aligned}
\]
The conclusion follows because the topologies of $\A_{Y}$ and $\AXp{\p}$ are characterized by commensurability with $\A_{Y}^{+}$ and $\AXp{\p}^{+}$ (\S~\ref{subsec:adeles} and \S\ref{subsubsec:the restricted direct product topology}).

\item $(3) \implies (1)$: Recall that $\Omega$ is discrete in $\A_{Y}$ (Remark~\ref{R:Sigma discrete in AX}). Thus, if the adelic topology on $\A_{Y}$ is the initial topology of $\phi : \A_{Y} \to \AXp{\p}$, there is some neighborhood $U_{D}$ of $0$ such that $\Omega \cap \phi^{-1}(U_{D}) = \{0\}$. Then $\phi(\Omega \cap \phi^{-1}(U_{D})) = \phi(\Omega) \cap U_{D} = \{0\}$. In view of Lemma~\ref{L:AY -> AXp continuity} this says that $\Omega$ is discrete in $\AXp{\p}$.
% The result is easily seen to hold in general for topological groups.
% 
\end{itemize}
This shows the equivalence of the three conditions. Now, assuming they hold, let $y \in \pi^{-1}(x)$. By Lemma~\ref{L:AY -> AXp injectivity}, we have   $y = s_j(x)$ for some (not necessarily unique) $j : 1 \leq j \leq m_{x}$. Thus we have a tower of field extensions $K_{x} \hookrightarrow L_{y} \hookrightarrow K_{x,j}$. Hence
\[
	     \deg_{\Sigma} \Omega
	=    \!\!\!\sum_{\pi(y)=x}\dim_{K_{x}}L_y
	\leq \sum_{j=1}^{m_{x}}\dim_{K_{x}}K_{x,j}
	=    \dim_{K_{x}} \Kxp{\p_{x}} 
	=    n,
\]
since $\p_{x}$ is monic of degree $n = \deg \p$.
\end{proof}

\begin{example}[A non-discrete subfield of an adelic algebra]
\label{EX:non-discrete subfield}	
Let $X=\mathbb{P}_1(\k)$ where $\chr(\k) \neq 2$, and $\Sigma=\k(u)$ with $u$ an unknown. Let $f\in \Sigma^{*} \setminus \Sigma^{*2}$, $Y$ the Zariski-Riemann variety of $\Sigma[T]/(T^2-f)$ (which is irreducible), and 
$\pi:Y\to X$ the induced map. Fix a point $x_{0} \in X$ such that $\pi^{-1}(x_{0}) = \{y_{1}, y_{2}\}$ consists of two distinct points $y_{j} = s_{j}(x_{0})$. Clearly the fields $K_{x_{0},1}$ and $K_{x_{0},2}$ are isomorphic. Fix an isomorphism $\varphi_{0}: K_{x_{0},2} \isomto K_{x_{0},1}$ of $K_{x_{0}}$-algebras. Finally, denote by $\f \in \I_{X}$ the idele of function germs $\f = (f_{x})_{x \in X}$ and $\p(T) = T^{2} - \f \in \A_{X}[T]$.

Let $\Omega$ be the image of the map $\phi: \Sigma[T]/(T^2-f) \to \AXp{\p} = \A_{X}[T]/(T^2-\f)$ defined by
\begin{equation}
	g \mapsto \phi(g) = (\phi(g)_y)_{y\in Y}
	\eqdef
	\begin{cases} 
		g_{y} & y \neq y_{1} \\
		\varphi_{0}(g_{y_{2}}) & y=y_{1},
	\end{cases}
\end{equation}
where $g_y$ denotes the germ of $g$ at $y\in Y$. Then, since the infinite dimensional vector space  $H^0(Y\setminus\{y_{1}\},\oo_{Y}) $ lies inside $\Omega \cap \AXp{\p}^{+}$, it follows that $\Omega$ is a subfield that is not discrete.
\end{example}

% (end)
% end subsection discrete subfields

%%%%%%%%%%%%%%%%%%%%%%%%%%%%%%%%%%%%%%%%%%%%%%%%%%%%%%%%%%%%%%%%%%%%%%%%%%%%%%%%%%%
%%%%%%%%%%%%%%%%%%%%%%%%%%%%%%%%%%%%%%%%%%%%%%%%%%%%%%%%%%%%%%%%%%%%%%%%%%%%%%%%%%%
\subsection{Product formula}
\label{subsec:product formula}
% (fold)
%%%%%%%%%%%%%%%%%%%%%%%%%%%%%%%%%%%%%%%%%%%%%%%%%%%%%%%%%%%%%%%%%%%%%%%%%%%%%%%%%%%
%%%%%%%%%%%%%%%%%%%%%%%%%%%%%%%%%%%%%%%%%%%%%%%%%%%%%%%%%%%%%%%%%%%%%%%%%%%%%%%%%%%

In attempting to axiomatize class field theory, Artin and Whaples~\cite{ArtinW} pointed out the importance of the product formula for valuations. They showed that a global field may be characterized by the existence of a product formula, along with having either an archimedean valuation or one with a finite residue field.

Here we use an additive analog of the usual notion of content, which is expressed multiplicatively in terms of valuations. This analog uses the definition of index (see~\cite[\S2.2.4]{AndersonPablos}) introduced when studying reciprocity laws for function fields of curves.

The main result is contained in Theorem~\ref{T:product formula equivalence}, which adds to the equivalent characterizations of Theorem~\ref{T:equivalences inyectivity} a new zero-content condition.

% %%%%%%%%%%%%%%%%%%%%%%%%%%%%%%%%%%%%%%%%%%%%%%%%%%%%%%%%%%%%%%%%%%%%%%%%%%%%%%%%%%%
% \subsubsection{Content function on an adelic algebra}
% \label{subsubsec:content}
% % (fold)
% %%%%%%%%%%%%%%%%%%%%%%%%%%%%%%%%%%%%%%%%%%%%%%%%%%%%%%%%%%%%%%%%%%%%%%%%%%%%%%%%%%%

\begin{definition}[Content function on the unit group of $\AXp{\p}$]
\label{D:degree invertible}
Given $\alpha \in \AXp{\p}^{*}$, we define its \emph{content} as
\begin{multline}
\label{E:degree invertible}
	\content(\alpha) \eqdef 
	\dim_{\k}\left(\faktor{\AXp{\p}^{+} }{\AXp{\p}^{+} \cap \alpha \AXp{\p}^{+}} \right)
	\\
	- \dim_{\k} \left(\faktor{\alpha \AXp{\p}^{+}}{ \AXp{\p}^{+} \cap \alpha \AXp{\p}^{+}} \right).
\end{multline}
Note that this only depends on the commensurability class of $\AXp{\p}^{+}$.
\end{definition}

That the definition makes sense follows from the fact that the adelic topology is an $\A_{X}$-algebra topology which can be characterized via commensurability (Lemma~\ref{L:conmensurability topology equivalence}). In particular multiplication by an element $\alpha \in \AXp{\p}$ is continuous, and thus multiplication by an invertible element is a homeomorphism. Hence the subspaces $\AXp{\p}^{+}$ and $\alpha \AXp{\p}^{+}$ are commensurable.

% \begin{corollary}
% \label{C:commensurabilityInvertibles}
% Let $\p(T) \in \A_{X}[T]$ be monic and separable. Given an invertible element $\alpha \in \AXp{\p}$, the subspaces $\AXp{\p}^{+}$ and $\alpha \AXp{\p}^{+}$ are commensurable.
% \end{corollary}

\begin{lemma}
\label{L:content computed pointwise}
For $\alpha \in \AXp{\p}^{*}$ with $\alpha = (\alpha_{x})$, we have
\begin{equation}
\label{E:content computed pointwise}
	\content(\alpha) =
	\sum_{x \in X} \left( \dim_{\k} \left(\faktor{\bb_{x}}{\bb_{x} \cap \alpha_{x} \bb_{x}} \right) 
	- \dim_{\k} \left(\faktor{\alpha_{x} \bb_{x}}{ \bb_{x} \cap \alpha_{x} \bb_{x}} \right) \right)
\end{equation}
where almost all summands are null. In terms of the decomposition~\eqref{E:KxKj BxBj} we have
\begin{equation}
\label{E:content sum valuations vxj}
	\content(\alpha) = \sum_{(x,j)} \upsilon_{x,j}(\alpha_{x,j})
\end{equation}
where $\alpha_{x,j}$ is the image of $\alpha_{x}$ in $ K_{x,j}$.
\end{lemma}

\begin{proof}
Since $\AXp{\p}^{+} = \prod_{x} \bb_{x}$ by definition, we have
\[
	\begin{aligned}
	   \faktor{\AXp{\p}^{+} }{\AXp{\p}^{+} \cap \alpha \AXp{\p}^{+}}
	&= \prod_{x} \left(\faktor{\bb_{x}}{\bb_{x} \cap \alpha_{x} \bb_{x}} \right)
	 % = \prod_{(x,j)} \left(\faktor{B_{x,j}}{B_{x,j} \cap \alpha_{x,j} B_{x,j}} \right)
	\\
	   \faktor{\alpha \AXp{\p}^{+} }{\AXp{\p}^{+} \cap \alpha \AXp{\p}^{+}}
	&= \prod_{x} \left(\faktor{\alpha_{x} \bb_{x}}{\bb_{x} \cap \alpha_{x} \bb_{x}} \right)
	 % = \prod_{(x,j)} \left(\faktor{\alpha_{x,j} B_{x,j}}{B_{x,j} \cap \alpha_{x,j} B_{x,j}} \right)
	\end{aligned}
\]
and the left hand sides are finite-dimensional $\k$-vector spaces since, as mentioned above, $\AXp{\p}^{+}$ and $\alpha \AXp{\p}^{+}$ are commensurable. Thus all factors on the right are also finite-dimensional, and indeed trivial for almost all $x$. Thus~\eqref{E:content computed pointwise} follows. Now, by~\eqref{E:KxKj BxBj},
\[
	\faktor{\bb_{x}}{\bb_{x} \cap \alpha_{x} \bb_{x}}
  =	\prod_{j=1}^{m_{x}} \faktor{B_{x,j}}{B_{x,j} \cap \alpha_{x,j} B_{x,j}},
  % \quad
  % 	\faktor{\alpha_{x} \bb_{x}}{\bb_{x} \cap \alpha_{x} \bb_{x}}
  % = \prod_{j=1}^{m_{x}} \faktor{\alpha_{x,j} B_{x,j}}{B_{x,j} \cap \alpha_{x,j} B_{x,j}}
\]
and similarly for ${\alpha_{x} \bb_{x}}/{\bb_{x} \cap \alpha_{x} \bb_{x}}$, from which~\eqref{E:content sum valuations vxj} follows.
\end{proof}

% \begin{multline}
% \label{E:content computed pointwise}
% 	\content(\alpha) =
% 	\sum_{x \in X} \deg(x) \cdot \left( \dim_{\k(x)} \left(\faktor{\bb_{x}}{\bb_{x} \cap \alpha_{x} \bb_{x}} \right) \right.
% 	\\
% 	\left.
% 	- \dim_{\k(x)} \left(\faktor{\alpha_{x} \bb_{x}}{ \bb_{x} \cap \alpha_{x} \bb_{x}} \right) \right),
% \end{multline}
% %%
% where $\deg(x) \eqdef \dim_{\k} \k(x)$.

In particular, multiplication by any $\alpha\in \I_{X}$ is a homeomorphism. For $n = 1$ and $\p(T) = T$ we have $\AXp{\p} \simeq \A_{X}$ and~\eqref{E:content computed pointwise} reduces to
\begin{equation}
\label{E:content over IX}
	\content(\alpha) = \sum_{x \in X} \upsilon_{x}(\alpha_{x})
	\qquad
	(\alpha \in \I_{X}),
\end{equation}
which is the additive analogue of the definition of content for the adeles of a global field (see \cite[Ch. II \S16]{CaFr}). Thus, the content of a function $f \in \Sigma^{*}$ is the degree of its divisor, in other words
\begin{equation}
\label{E:V=0 on Sigma}
	\content(f) = 0 \qquad (f \in \Sigma^{*}).
\end{equation}
As one might expect, it can also be shown (see~\cite[\S2.2.4]{AndersonPablos} for this in a more general context) that the content function is additive, namely
\begin{equation}
\label{E:}
	\content(\alpha\beta) = \content(\alpha) + \content(\beta)
	\qquad
	(\alpha, \beta \in \AXp{\p}^{*}).
\end{equation}

\begin{theorem}
\label{T:product formula equivalence}
Let $\Omega$ be a finite separable field extension of $\Sigma$ contained in  $\AXp{\p}$. Then the conditions in Theorem~\ref{T:equivalences inyectivity} are equivalent to the ``product formula''
\begin{equation}
\label{E:product formula}
	\text{$\content(f) = 0$ for all $f \in \Omega^{*}$.}
\end{equation}
\end{theorem}

\begin{proof}
Suppose that the map $\phi : \A_{Y} \to \AXp{\p}$ from~\eqref{E:AY -> AXp} is injective. Then~\eqref{E:injectivity Bx prod Ly = prod By} holds. To simplify notation, we will identify the corresponding objects with their images under $\phi$ and $\phi_{x}$. Thus, for $f \in \Omega^{*}$ and $x \in X$, taking products over $y \in \pi^{-1}(x)$, we have
\begin{align*}
	&\mathbin{\phantom{=}}
	   \faktor{\bb_{x} \cap \textstyle{\prod} L_{y}}{f_{x} \bb_{x} \cap \bb_{x} \cap \textstyle{\prod} L_{y}}
	= \faktor{\textstyle{\prod} B_{y}}{f_{x} \bb_{x} \cap \bb_{x} \cap \textstyle{\prod} L_{y} \cap f_{x} \textstyle{\prod} L_{y}}
	 % = \faktor{\textstyle{\prod} B_{y}}{f_{x} \bb_{x} \cap \textstyle{\prod} B_{y}}
	\\[1ex]
	&= \faktor{\textstyle{\prod} B_{y}}{f_{x} (\bb_{x} \cap \textstyle{\prod} L_{y}) \cap (\bb_{x} \cap \textstyle{\prod} L_{y})}
	 = \faktor{\textstyle{\prod} B_{y}}{f_{x} (\textstyle{\prod} B_{y}) \cap (\textstyle{\prod} B_{y})}
	% &= \faktor{\textstyle{\prod} B_{y}}{f_{x} (\prod B_{y})\bb_{x} \cap \textstyle{\prod} B_{y}}
	%  = \faktor{\textstyle{\prod} B_{y}}{\textstyle{\prod} f_{y} B_{y} \cap \textstyle{\prod} B_{y}}
	\\[1ex]
	&= \prod \faktor{B_{y}}{f_{y} B_{y} \cap B_{y}}.
\end{align*}
The first equality follows from observing that $f_{x} \prod L_{y} = \prod L_{y}$ since $f \in \Omega^{*}$. In general, we have the exact sequence
\begin{multline*}
	  0
	  \longrightarrow \faktor{\bb_{x} \cap \textstyle{\prod} L_{y}}{f_{x} \bb_{x} \cap \bb_{x} \cap \textstyle{\prod} L_{y}}
	\\[1ex]
    \longrightarrow \faktor{\bb_{x}}{\bb_{x} \cap f_{x} \bb_{x}}
	\longrightarrow \faktor{\bb_{x}}{\bb_{x} \cap \textstyle{\prod} L_{y}}
	\longrightarrow	0,
\end{multline*}
which by the previous reasoning is the same as
\begin{multline*}
	  0 
	  \longrightarrow \prod \faktor{B_{y}}{f_{y} B_{y} \cap B_{y}}
	  \longrightarrow \faktor{\bb_{x}}{\bb_{x} \cap f_{x} \bb_{x}}
	  \longrightarrow \faktor{\bb_{x}}{\textstyle{\prod} B_{y}}
	  \longrightarrow 0.
\end{multline*}
A similar argument leads to the second exact sequence
\begin{multline*}
	  0 
	  \longrightarrow \prod \faktor{f_{y} B_{y}}{f_{y} B_{y} \cap B_{y}}
	  \longrightarrow \faktor{f_{x} \bb_{x}}{\bb_{x} \cap f_{x} \bb_{x}}
	  \longrightarrow \faktor{f_{x} \bb_{x}}{\textstyle{\prod} f_{y} B_{y}}
	  \longrightarrow 0.
\end{multline*}
Note that all terms in the previous two sequences are finite-dimensional $\k$-vector spaces, and the last terms in each are isomorphic. Thus, summing over $x \in X$ and recalling~\eqref{E:content computed pointwise} and~\eqref{E:V=0 on Sigma}, we conclude that the content of $f$ considered as an element of $\AXp{\p}^{*}$ is the same as its content as an element of $\A_{Y}^{*} = \I_{Y}$, namely $0$.

For the converse, suppose that the map~\eqref{E:AY -> AXp} is not injective. Recalling the decomposition~\eqref{E:KxKj BxBj} and Lemma~\ref{L:AY -> AXp injectivity}, there is then some point $y \in Y$ such that no $\upsilon_{x,j}$ restricts to $\upsilon_{y}$, where $x = \pi(y)$. There is a nonconstant function $f$ which is regular in $Y \setminus \{y\}$ and has a zero at some $y_{0} = s_{j_{0}}(x_{0})$. Thus $\upsilon_{y}(f) < 0$, whereas $\upsilon_{x,j}(f) \geq 0$ for all $(x,j)$ and $\upsilon_{x_{0}, j_{0}}(f) > 0$, hence by~\eqref{E:content sum valuations vxj} it is clear that $\content(f) > 0$.
\end{proof}

% (end)
% subsection product formula end

\begin{remark}
\label{R:algebras instead of fields}
The reader may have noticed that the methods outlined in the previous sections may be readily generalized to the case when the field $\Omega$ is replaced by a finite $\Sigma$-algebra contained in an adelic algebra $\AXp{\p}$, and $\Omega^{*}$ is the group of invertible elements. Indeed, these are isomorphic to finite products of field extensions of $\Sigma$ and thus may be dealt with by reducing to that case. Theorem~\ref{T:equivalences inyectivity} and Theorem~\ref{T:product formula equivalence} hold in this more general setting.
\end{remark}

% (end)
% end section subfields of adelic algebras

%%%%%%%%%%%%%%%%%%%%%%%%%%%%%%%%%%%%%%%%%%%%%%%%%%%%%%%%%%%%%%%%%%%%%%%%%%%%%%%%%%%
%%%%%%%%%%%%%%%%%%%%%%%%%%%%%%%%%%%%%%%%%%%%%%%%%%%%%%%%%%%%%%%%%%%%%%%%%%%%%%%%%%%
%%%%%%%%%%%%%%%%%%%%%%%%%%%%%%%%%%%%%%%%%%%%%%%%%%%%%%%%%%%%%%%%%%%%%%%%%%%%%%%%%%%
\section{Concluding remarks}
\label{sec:concluding remarks}
% (fold)
%%%%%%%%%%%%%%%%%%%%%%%%%%%%%%%%%%%%%%%%%%%%%%%%%%%%%%%%%%%%%%%%%%%%%%%%%%%%%%%%%%%
%%%%%%%%%%%%%%%%%%%%%%%%%%%%%%%%%%%%%%%%%%%%%%%%%%%%%%%%%%%%%%%%%%%%%%%%%%%%%%%%%%%
%%%%%%%%%%%%%%%%%%%%%%%%%%%%%%%%%%%%%%%%%%%%%%%%%%%%%%%%%%%%%%%%%%%%%%%%%%%%%%%%%%%

As we mentioned in the introduction, for a cover of algebraic curves $Y \to X$, the primitive element theorem for the corresponding extension of their function fields $\Omega/\Sigma$ implies that $\A_{Y} \simeq \A_{X} \otimes_{\Sigma} \Omega$ is also generated by a primitive element. Clearly, an adelic algebra $\AXp{\p} = \A_{X}[T]/(\p(T))$ has the class of $T$ as a primitive element. However, there will also be others, and in particular, regarding the problem of determining when such an algebra arises from a cover, one is interested in primitive elements with minimal polynomial having coefficients not just in $\A_{X}$ but in the function field $\Sigma$. We will take up this problem in future work based on the machinery developed here, e.g. Proposition~\ref{P:Hensel complete splitting}, which may be used to give more explicit criteria characterizing such primitive elements. The existence of primitive elements in ring extensions is a long-standing question which has been studied, for example, in~\cite{AndersonDobbsMullinsC} and~\cite{BagioPaques}, although adele rings present technical difficulties, such as having infinitely many idempotents.

Another avenue of research is the study of reciprocity laws on algebraic curves, in the spirit of~\cite{AndersonPablos} or~\cite{MP} and generalizing the results of~\cite{MPNPP-collectanea}. In this context, Theorem~\ref{T:Axp is free AX-module}, which shows that $\AXp{\p}$ is a free topological $\A_{X}$-module, leads one to consider the Grassmannian $\operatorname{Gr}(\A_{X}^{\oplus n})$. A group action on 
$\operatorname{Gr}(\A_{X}^{\oplus n})$ gives a central extension of the group. It should be straightforward to check that discreteness of a function field $\Omega$ lying in $\AXp{\p}$ and the product formula determine a point in the Grassmannian fixed under the action of $\Omega^{*}$ and implies the triviality of the associated central extension, which yields a global symbol and expresses a corresponding reciprocity law. It should also be interesting to study how this varies for different extensions $\Omega/\Sigma$.

% (end)
% end section concluding remarks

\end{document}